# Split-Step Method for Generalized Non-Linear Equations: A Three Operator Story

Haider Zia[1]*

[1]*Max-Planck Institute for the Structure and Dynamics of Matter, Luruper Chausee 149, Hamburg, DE-22607*
*Departments of Chemistry and Physics, University of Toronto*

*\*haider.zia@mpsd.mpg.de*

**Abstract:** This paper describes an updated exponential Fourier based split-step method that can be applied to a greater class of partial differential equations than previous methods would allow. These equations arise in physics and engineering, a notable example being the generalized derivative non-linear Schrödinger equation that arises in non-linear optics with self-steepening terms. These differential equations feature terms that were previously inaccessible to model accurately with low computational resources.  The new method maintains a $3^{rd}$ order error even with these additional terms and models the equation in all three spatial dimensions and time. The class of non-linear differential equations that this method applies to is shown. The method is fully derived and implementation of the method in the split-step architecture is shown. This paper lays the mathematical ground work for an upcoming paper employing this method in white-light generation simulations in bulk material.

## 1 Niche of the New General Split-Step Method

This paper presents an updated methodology that is based on the Strang exponential Fourier split-step method (EFSSM) [1]  that can model complicated non-linear partial differential equations (PDE) in all three spatial dimensions and time. These non-linear equations contain a distributional coefficient of a derivative term (the distributional coefficient is over the same variable(s) as the derivative term). This paper provides the mathematical support of a new white light generation simulation which is presented in [2].

Amongst a wide plethora of fields, these non-linear PDEs mostly arise in physical sciences and engineering disciplines. A notable example includes the generalized derivative non-linear Schrödinger equation in non-linear optics [3,4] that incorporates self-steepening effects. This equation will be shown in the text as an example equation for the application of the method. More examples are other derivative non-linear Schrödinger equations in plasma physics [5], fluid dynamics (Benjamin-Bona-Mahony equation, and water waves [6,7]), integrable systems and dynamical Systems.

Current spectral split-step methods have been extensively explored in one dimensional systems such as the one-dimensional cubic non-linear Schrödinger Equation (NLSE). Past studies focused primarily on estimating step-size dependent error and deriving adaptive step-size algorithms for the implementation of the method [8-10]. Other studies explored the stability of the method [11] or its application around certain bound solutions of a NLSE such as soliton formation [12]. However, extensions of these methods have not been comprehensively explored until now. Exponential Fourier split-step methods need to be applicable to beyond cubic NLSE type equations and in all spatial dimensions and time. Because, the exponential nature of this approach reduces the amount of operations necessary and increases accuracy as compared to finite difference and Runge Kutta methods. Since, it more closely matches the "analytic" integrated solutions of the corresponding non-linear PDE. The computational time, the implementation ease, the low amount of operations and the accuracy of the newly derived extended EFSSM that will be presented offers an attractive alternate to other numerical schemes.

In summary, a new methodology that can accurately simulate the additional complexity of non-linear PDEs that will be introduced below, in all spatial dimensions and time, will be derived. The exponential split-step Fourier approach is taken as a basis and a new extension is derived.

## 2 Derivation of the General Method

Split-step methods are based on modelling the non-linear PDE with operators acting on the solution in an iterative fashion in the propagation coordinate. Terms of the non-linear PDE become operators. Per propagation step, the non-linear equation is decomposed in a series of coupled equations in terms of each operator (named flows of the equation). These equations are then numerically or analytically solved. The method scales with 3rd order accuracy in the propagation coordinate.

To start the derivation of the new extended method, a discussion of the type of solutions that the method aims to find will be carried out and the general form of non-linear partial differential equations the method can be applied to will be presented. Operators will be defined and the application of these operators will be discussed, especially the new operators introduced in this extended method. The Fourier domains in which these operators are applied in for the best numerical accuracy will be shown.

### 2.1 Non-linear Partial Differential Equation Types Supported by Method

The method models integrable bound L-2 normed [13] solutions of the corresponding non-linear PDE. The necessity of this condition is because the Fourier integral for these functions always exists, which is needed since this method uses the Fourier transform and its properties. The general form of non-linear equations of which is considered writes as:

$$\frac{\partial u}{\partial \varsigma} = \wp u + \mathbb{Q}\mathbb{N}u \qquad (1)$$

This equation is defined in a unit-less coordinate system and for unit-less quantities. Where, $x, y$ are the transverse unit-less coordinates, $\tau$ is the unit-less time coordinate. $\varsigma$, is the unit-less propagation coordinate.

$\wp$ is composed of a power series over derivative operators $(\nabla, \frac{\partial}{\partial \tau})$, shown in Eq. (2). If the series of Eq. (2) converges to a closed form expression when the derivative operators are replaced by complex variables, $\wp$ can be described as this closed form function where the complex variable is replaced by the derivative operators. This is in accordance with the definitions used in operator analysis theory [14]. Appendix 0.2 discusses in more detail the consequences of when the series does not converge to a closed form functional form within the context of the complete methodology.

$$\wp = \sum_{n=0}^{\infty} \sum_{j=0}^{\infty} c_{nj} \nabla^n \left(\frac{\partial}{\partial \tau}\right)^j \qquad (2)$$

The operator $\nabla^n$ is defined as such: $\nabla^n \equiv \frac{\partial^n}{\partial^n x} + \frac{\partial^n}{\partial^n y}$ [1]. While, for the derivative operator $\frac{\partial}{\partial \tau}$: $\left(\frac{\partial}{\partial \tau}\right)^j = \frac{\partial^j}{\partial^j \tau}$ follows from the properties of the derivative. The $c_{nj}$ are complex constants.

$$\mathbb{Q} = \left(c_1 + c_2 \frac{\partial}{\partial \tau}\right) \quad (3)$$

$c_1, c_2$ are complex constants.

$$\mathbb{N} = \beta(u, v(x_\perp, \tau)) \quad (4)$$

β is a function whose arguments are the independent coordinate variables in space and time and distributions, e.g., based on the absolute value of $u$ to an arbitrary $j^{th}$ order. Where, the notation $x_\perp$ to mark the set of unit-less transverse coordinates is used. $v(x_\perp, \tau)$ represents a set of functions over $x_\perp, \tau$.

The reason for why such a specific form of non-linear partial differential equations is considered here is that this corresponds to pertinent generalized derivative non-linear Schrödinger equations central to the study of non-linear optics. For example, the methodology developed in this paper was used to simulate the non-linear propagation of light in bulk material described in [2], using the equation factoring all pertinent higher order effects described by [3]. The equation used in [2] is of the general form:

$$\frac{\partial u}{\partial \varsigma} = \frac{i}{4}\left(1 + \frac{i}{a}\frac{\partial}{\partial \tau}\right)^{-1} \nabla_\perp^2 u - ib \frac{\partial^2 u}{\partial \tau^2} + i\left(1 + \frac{i}{a}\frac{\partial}{\partial \tau}\right)\left[c|u|^2 u - d\left(1 - \frac{i}{e}\right)\rho u + if|u|^{2(m-1)}u\right] \quad (5)$$

$a - f$ are physical constants and $\rho$ is a function based on $u^2$. The first two terms, $\frac{i}{4}\left(1 + \frac{i}{a}\frac{\partial}{\partial \tau}\right)^{-1} \nabla_\perp^2 u - ib \frac{\partial^2 u}{\partial \tau^2}$, can be described in the following series representation[3], and is, therefore, in the form of $\wp$ shown in Eq. (2):

$$\frac{i}{4}\left(1 + \frac{i}{a}\frac{\partial}{\partial \tau}\right)^{-1} \nabla_\perp^2 u - ib \frac{\partial^2 u}{\partial \tau^2} = \sum_{n=0}^{\infty}\sum_{j=0}^{\infty}[c_n(-1)^j(\frac{i}{a})^j + d_{nj}](\frac{\partial}{\partial \tau})^j \nabla_\perp^n = \wp \quad (6)$$

Where, $c_n = 0$ when $n \neq 2$, otherwise $c_n = \frac{i}{4}$. $d_{nj} = 0$ if $n \neq 0$ or $j \neq 2$, otherwise, $d_{nj} = \frac{-ib}{2}$.

As well, since $\left[c|u|^2 u - d\left(1 - \frac{i}{e}\right)\rho u + if|u|^{2(m-1)}u\right]$ is in the form of $\mathbb{N}$ and $i\left(1 + \frac{i}{a}\frac{\partial}{\partial \tau}\right)$ is in the form of $\mathbb{Q}$. Therefore, Eq. (5) is of the same form as Eq. (1).

---

[1] $\frac{\partial^0}{\partial^0 x}, \frac{\partial^0}{\partial^0 y}, \frac{\partial^0}{\partial^0 \tau} \equiv 1$. Also, the defining series for $\wp$ can be rewritten over a triple summation separating the x and y derivative terms. This definition would cover more cases and the method presented in this paper can equally be applied to this more general definition. However, to reduce the amount of summations (for ease of legibility), while still fully covering the case presented in [2] this definition is used.

[2] The derivation of Eq. (5) and the physics involved is explained at great length in [3], [2]. Here the equation is presented to demonstrate why the form of Eq. (1) is important to consider.

[3] Using the complex valued binomial expansion of $\left(1 + \frac{i}{a}\frac{\partial}{\partial \tau}\right)^{-1}$.

There are two reasons for why $\wp$ is defined as in Eq. (2) and not in a closed form function. The first reason being, to increase the breadth of relevance to other useful generalized non-linear Schrödinger equations. The second reason will become clear later in this section.

## 2.2 Operator Definitions and Representation in Relevant Domains

The coefficients of Terms in Eq. (1) will now be grouped into a series of operators that are applied in an iterative sense along the propagation coordinate in the goal of producing a solution taking $u$ inputted into the first slice of the propagation coordinate as the initial and boundary conditions of the equation. In order to accomplish this goal, the operators will first be defined, then a method for how they are applied will be derived with 3$^{rd}$ order accuracy or higher. The operators are chosen such that they are accurate and simple in their application. It is shown later (section 2.3) that this is accomplished if terms in Eq. (1) are grouped together into operators by the domains of application. Since the structure of the term depends on its domains of application, terms in an operator happen to have the same structure as well. For example, coefficient terms of Eq. (1) that are composed of a constant coefficient multiplied into a differential operator are grouped together into an operator because these terms are simplified in the Fourier domains (i.e., by domains of applicability). Other coefficient terms that are composed only of functions and distributions over complex numbers (and not operators) are grouped together into another operator because they are all applied in the original domains. The third operator consists of the remaining terms which are composed of a functional or distributional coefficient over complex numbers multiplied into a derivative operator. It will be shown that these terms are applied in a mixture of Fourier and original domains and are of a more complicated nature than the previous two operator types. This third operator will be explored extensively in section 2.3-2.5.

The coefficient $\wp$ is taken as the straight forward linear operator. This corresponds to the first operator described in the previous paragraph. The terms in $\wp$ are grouped together into one operator because it can be shown that the derivative operators can be replaced by imaginary numbers if this operator is applied on $u$ in the Fourier space. Also, this is true for the exponential version of this operator which is used to model the terms of Eq. (1) whose coefficients are covered within this operator (the way it is applied is shown rigorously in section 2.3). This is of great benefit since in its Fourier space $\wp$ is a function over numerical variables instead of derivative operators and can be easily evaluated. This is demonstrated in Appendix 0.1-0.3. This is also in accordance with the use of the Fourier representation of the linear operator as is done in traditional pseudo-spectral techniques [9]. Therefore, -i$w$ is substituted for $\frac{\partial}{\partial \tau}$ and $-ik_x, -ik_y$ for $\nabla_x^1, \nabla_y^1$ (($k_x, k_y, w$) are the angular frequencies of the time and spatial variables) obtaining $\wp$ in the inverse space:

$$\wp(k_x k_y, w) = \sum_{n=0}^{\infty} \sum_{j=0}^{\infty} c_{nj}((-ik_x)^n + (-ik_y)^n)(-iw)^j \qquad (7)$$

Where, the above series can converge to a closed form functional representation over a certain range in $(k_x, k_y, w)$. Appendix 0.2 lists convergence conditions for when it is possible to have a closed form functional representation in the inverse domain. This easy conversion of derivative operators to complex number variables is the second reason for why $\wp$ was defined in its series representation from section 2.1.

$\mathbb{QN}$ is split into distinct operators labelled as $\alpha_1, \alpha_2$. $\mathbb{QN}$ can be expanded as:

$$\mathbb{Q}\mathbb{N} = \left(c_1\mathbb{N} + c_2\left(\frac{\partial \mathbb{N}}{\partial \tau} + \mathbb{N}\frac{\partial}{\partial \tau}\right)\right) = \alpha_1 + \alpha_2 \tag{8}$$

$$\alpha_1 = c_1\mathbb{N} + c_2\left(\frac{\partial \mathbb{N}}{\partial \tau}\right) \tag{9}$$

And,

$$\alpha_2 = c_2\mathbb{N}\frac{\partial}{\partial \tau} \tag{10}$$

$\mathbb{Q}\mathbb{N}$ is decoupled into an operator solely defined on distributions, and an operator defined by a distribution and a derivative operator. This decoupling provides an accurate way of modelling the total term as will be seen later in this section. The $\alpha_2$ operator is novel and the operator and the method of its application is not present in previous split-step methods.

$\alpha_1$ is a non-linear operator and composed of functions and distributions over complex numbers and not on derivative operators. $\alpha_2$ is named the augmented non-linear operator. It is composed of a distributional coefficient over complex numbers to a derivative operator.

In contrast to $\wp$, it would be of no benefit to consider $\alpha_1$ in any inverse space, since it is already over complex numerical variables and not operators. Hence, its exponential form used in section 2.3 is applied in the original space $\alpha_1$ was defined in. If $\alpha_1$ is defined over derivatives of functions of $u$, as is the case if using Eq. (5) and in [2], the derivatives are evaluated in an analytic way with Fourier operations; The Fourier identity for the derivative outlined above is used and the term is transformed back into the time domain. For functions/distributions in $\alpha_1$ that are defined in terms of satisfying an additional differential equation, as is also the case in [2], the appropriate time domain integration method is used, i.e. Runge-Kutta, etc. It will later be demonstrated that the complexity of $\alpha_2$ can be reduced and then suitable domains to apply it in will be found, as done for $\wp, \alpha_1$.

Now, that operators corresponding to terms in Eq. (1) are defined, the next section will derive the method in how they are applied to simulate Eq. (1) with third order accuracy.

## 2.3 Flows of the Non-linear Partial Differential Equation under the 3 Operator Decomposition

This section will now derive and show how a given chosen set of operators (that cover all terms in Eq. (1), in this case chosen as in what is described in section 2.2) are applied within their domains of application to solve Eq. (1). Eq. (1) is broken up into a series of differential equations involving $\wp, \alpha_1, \alpha_2$. Meaning that each operator (consisting of coefficient terms in Eq. (1)) acts as the coefficient of $u$ in its own first order ordinary differential equation w.r.t the propagation coordinate. Each equation is the named the 'flow' of the operator. Analytic or near analytic partial solutions $u_\wp, u_{\alpha_1}, u_{\alpha_2}$, for these differential equations involving $\wp, \alpha_1, \alpha_2$ are obtained. These equations are solved in a specific order. The solution of one equation acts as the initial value (in terms of the propagation coordinate) of the next equation. When applied in the shown specific order, within the step being considered, they yield an approximation of $u$ coming out of the propagation slice into the next slice. This methodology was adapted from [8, Eq.8-10] to account for three operators instead of two. The specific ordering, named Strang symmetrisation [1] is needed because each of the three operators have a nonzero commutation relation to the others. Therefore, their ordering influences the result introducing numerical error. The step-size error is then proportional to the square of the propagation step size [1]. The flows are now shown in the specific ordering they are applied in [2]:

1. $\begin{cases} \frac{\partial u_\wp(\varsigma)}{\partial \varsigma} = \wp u_\wp(\varsigma), \ \forall \varsigma \in [\varsigma_k, \varsigma_{k+\frac{1}{4}}] \\ u_\wp(\varsigma_k) = u(\varsigma_k), \ u \ coming \ from \ previous \ slice \end{cases}$

2. $\begin{cases} \frac{\partial u_{\alpha_2}(\varsigma)}{\partial \varsigma} = \alpha_2 u_{\alpha_2}(\varsigma), \ \forall \varsigma \in [\varsigma_k, \varsigma_{k+\frac{1}{4}}] \\ u_{\alpha_2}(\varsigma_k) = u_\wp\left(\varsigma_{k+\frac{1}{4}}\right) \end{cases}$

3. $\begin{cases} \frac{\partial u_\wp(\varsigma)}{\partial \varsigma} = \wp u_\wp(\varsigma), \ \forall \varsigma \in [\varsigma_{k+\frac{1}{4}}, \varsigma_{k+\frac{1}{2}}] \\ u_\wp\left(\varsigma_{k+\frac{1}{4}}\right) = u_{\alpha_2}\left(\varsigma_{k+\frac{1}{2}}\right), \end{cases}$

4. $\begin{cases} \frac{\partial u_{\alpha_1}(\varsigma)}{\partial \varsigma} = \alpha_1 u_{\alpha_1}(\varsigma), \ \forall \varsigma \in [\varsigma_k, \varsigma_{k+1}] \\ u_{\alpha_1}(\varsigma_k) = u_\wp\left(\varsigma_{k+\frac{1}{2}}\right) \end{cases}$

5. $\begin{cases} \frac{\partial u_\wp(\varsigma)}{\partial \varsigma} = \wp u_\wp(\varsigma), \ \forall \varsigma \in [\varsigma_{k+\frac{1}{2}}, \varsigma_{k+\frac{3}{4}}] \\ u_\wp\left(\varsigma_{k+\frac{1}{2}}\right) = u_{\alpha_1}(\varsigma_{k+1}), \end{cases}$

6. $\begin{cases} \frac{\partial u_{\alpha_2}(\varsigma)}{\partial \varsigma} = \alpha_2 u_{\alpha_2}(\varsigma), \ \forall \varsigma \in [\varsigma_{k+\frac{1}{2}}, \varsigma_{k+1}] \\ u_{\alpha_2}\left(\varsigma_{k+\frac{1}{2}}\right) = u_\wp\left(\varsigma_{k+\frac{3}{4}}\right) \end{cases}$

7. $\begin{cases} \frac{\partial u_\wp(\varsigma)}{\partial \varsigma} = \wp u_\wp(\varsigma), \ \forall \varsigma \in [\varsigma_{k+\frac{3}{4}}, \varsigma_{k+1}] \\ u_\wp\left(\varsigma_{k+\frac{3}{4}}\right) = u_{\alpha_2}(\varsigma_{k+1}) \end{cases}$

Where the final $u_\wp(\varsigma_{k+1})$ calculated from the 7th step becomes $u$ entering the next slice labelled $(\varsigma_{k+1})$ in the propagation direction and the above 7 steps are repeated iteratively. Fig.1 visually represents the above iteration. The methodology of generating u, how the operators $\wp, \alpha_1, \alpha_2$ are applied, the ordering of how they are applied and to what they are applied to have just been shown. However, the above procedure can be simplified because the above differential equations can be solved analytically or "close to" analytically for all operators, yielding:

For the first step:

$$u_\wp(\varsigma) = e^{\wp(\varsigma - \varsigma_k)} u(\varsigma_k) \tag{11}$$

Coming out of the first step then:

$$u_\wp\left(\varsigma_{k+\frac{1}{4}}\right) = e^{\frac{1}{4}\wp\Delta\varsigma} u(\varsigma_k) \tag{12}$$

For the second step:

$$u_{\alpha_2}(\varsigma) = e^{\int_{\varsigma_k}^{\varsigma} \alpha_2(\varsigma')d\varsigma'} u_\wp\left(\varsigma_{k+\frac{1}{4}}\right) \tag{13}$$

The integral arises since the $\alpha_2$ operator contains functions/distributions over $\varsigma$. However, since the interval $[\varsigma_k, \varsigma_{k+\frac{1}{2}}]$ is considered to be small such that these functions/distributions are slowly varying. Therefore, a mean value theory [15] can be employed. Where, the value of $\alpha_2$ in the exponent is calculated with $u$ given from the end of the previous operator step is used.

This gives a value at $\varsigma_{k+\frac{1}{2}}$ (at the end of the step) of:

$$u_{\alpha_2}\left(\varsigma_{k+\frac{1}{2}}\right) = e^{\frac{1}{2}\alpha_2\left(u_\wp\left(\varsigma_{k+\frac{1}{4}}\right)\right)\Delta\varsigma} u_\wp\left(\varsigma_{k+\frac{1}{4}}\right) \tag{14}$$

At the end of the third step:

$$u_\wp\left(\varsigma_{k+\frac{1}{2}}\right) = e^{\frac{1}{4}\wp\Delta\varsigma} u_{\alpha_2}\left(\varsigma_{k+\frac{1}{2}}\right) \tag{15}$$

At the end of the fourth step, again employing a mean value theory and using the value of $u$ coming from the previous step for all functions of $u$ in $\alpha_1$ at $\varsigma_{k+1}$, the following is obtained:

$$u_{\alpha_1}(\varsigma_{k+1}) = e^{\alpha_1\left(u_\wp\left(\varsigma_{k+\frac{1}{2}}\right)\right)\Delta\varsigma} u_\wp\left(\varsigma_{k+\frac{1}{2}}\right) \tag{16}$$

At the end of the fifth step the following is obtained:

$$u_\wp\left(\varsigma_{k+\frac{3}{4}}\right) = e^{\frac{1}{4}\wp\Delta\varsigma}u_{\alpha_1}(\varsigma_{k+1}) \tag{17}$$

At the end of the sixth step:

$$u_{\alpha_2}(\varsigma_{k+1}) = e^{\frac{1}{2}\alpha_2\left(u_\wp\left(\varsigma_{k+\frac{3}{4}}\right)\right)\Delta\varsigma}u_\wp\left(\varsigma_{k+\frac{3}{4}}\right) \tag{18}$$

And at the end of the final step:

$$u(\varsigma_{k+1}) = u_\wp(\varsigma_{k+1}) = e^{\frac{1}{4}\wp\Delta\varsigma}u_{\alpha_2}(\varsigma_{k+1}) \tag{19}$$

The proof of the integration of these steps is discussed in appendices 0.4, 1.1, 1.2. The proof is constructed by showing that the Maclaurin series of the exponential functions shown above satisfy the corresponding differential equation. Once this is shown, the fact that this series converges to the above exponential forms concludes the proof. Appendix 1.1 additionally applies the proof of these steps to derive the exponential form of a convolution type operator.

The following scheme is equivalent to:

$$u(\varsigma_{k+1}) = e^{\frac{1}{4}\wp\Delta\varsigma}e^{\frac{1}{2}\alpha_2\Delta\varsigma}e^{\frac{1}{4}\wp\Delta\varsigma}e^{\alpha_1\Delta\varsigma}e^{\frac{1}{4}\wp\Delta\varsigma}e^{\frac{1}{2}\alpha_2\Delta\varsigma}e^{\frac{1}{4}\wp\Delta\varsigma}u(\varsigma_k) \tag{20}$$

Where, the non-linear operators that are dependent on functions of $u$ are calculated at $u$ outputted from the previous operator step.

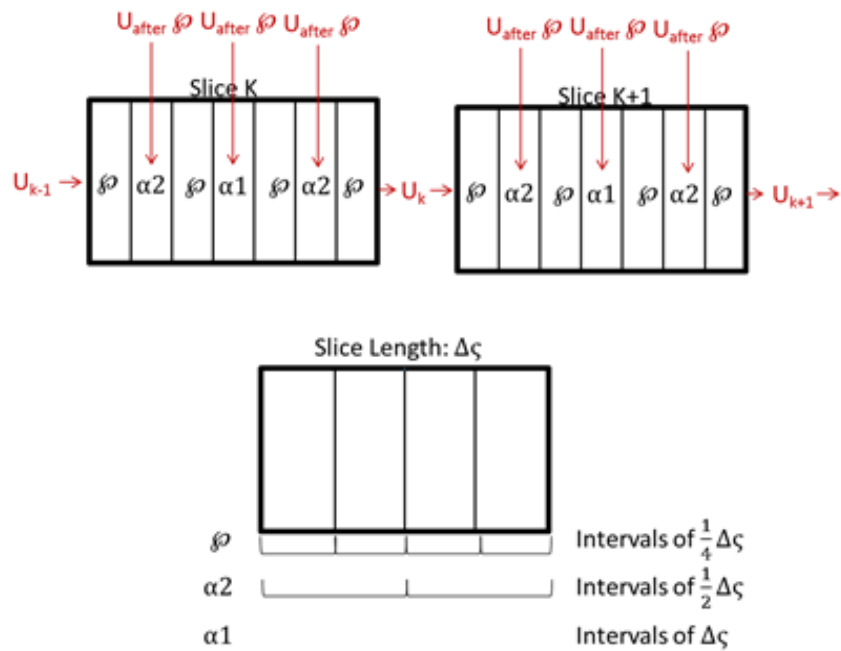

Figure 1: Graphical representation of steps 1-7. $u_{k-1}$ is calculated from the output of the slice numbered $k-1$ and is used as the initial $u$ in slice $K$. The updated $u$ from the previous iterative sub steps in slice $K$ is used to calculate the operator values for the pertinent sub step. Each operator is applied in intervals equivalent to the slice length divided by the amount of times the operator is applied in the slice.

The above ordering can be modified to fulfill different symmetrisations. In general, terms that are fast varying relative to others should be the ones most extensively split in the symmetrisation, as one would need to "sample $u$ interacting with these terms" more often. The symmetrisation is not a hard constraint and is dependent on the problem and the computational resources. Different symmetrisations have different convergence speeds with respect to the propagation coordinate step size.

From the above, the operator is an argument of an exponent. This is defined as the exponential form of the operator. The exponential form of each operator is multiplied with u. This is equivalent to the Maclaurin series expansion of the exponential with respect to the propagation coordinate (valid everywhere in that domain) multiplied into u. Appendix 0.3-0.4 outlines in more detail how this procedure is done. Appendix 1 details the proof for why the above differential equations integrate to the Maclaurin series of an exponent with the operator as an argument.

In the case where the argument operator has a range in the complex numbers, the terms of the Maclaurin series do not depend on u to be defined and can be determined (here as numbers) before they are applied to u. Also, since the argument operator equates to complex numbers over the domain variables, the exponent can be calculated without the need to use the Maclaurin series and it can be determined before being applied to **u**. I.e., the numerical value of the exponential can be computed directly and then multiplied into $u$ at the appropriate domain coordinates. Thus, the Maclaurin series does not need to be used in the application of this operator to u. However, in the case where the argument operator does not have a range in the complex numbers (for example, is composed of derivative operators), then each term of the Maclaurin series depends on u to be determined and the full Maclaurin series will have to be used and the exponential operator cannot be determined beforehand. Given the above discussion, it is prudent to find domains where the operators of the exponential arguments have a range in the complex numbers.

If the operator used in the exponential argument can be represented in a closed form expression, then the exponential argument can be computed with a finite amount of operations and no series truncation would have to occur. Therefore, it is prudent to work with operators that are closed form expressions with ranges over complex numbers within a certain domain, whenever possible.

It can be derived (contained in Appendix 0.1-0.3) that the Maclaurin series expansion of the exponential form of an operator that is a series of derivative terms with constant coefficients, multiplied into $u$ in the original space is equivalent to the inverse Fourier transform of the exponential of the representation of that operator in the Fourier space (i.e., frequency space) applied to the representation of $u$ in the Fourier space. Therefore, the exponential form of the $\wp$ operator in Fourier space will be the exponent of a function over independent numerical domain variables (instead of differential operators) and can be computed directly in this space.

Given this, it is not prudent to evaluate the Maclaurin series expansion of the exponential operator of $\wp$ in the original space. Because, each term of the Maclaurin series will consist of a series expansion of derivative operators acting on u. Therefore, numerical differentiation will have to be employed and both truncation in the terms of the Macluarin series and the series representation of each term will have to be employed.

If the representation of $\wp$ reduces to a closed form functional form, then the exponent can be directly calculated without need to calculate each term of the Maclaurin expansion (as discussed above). If not, each term in the Maclaurin expansion will have to be calculated in this Fourier space for each domain value and both truncation in the Maclaurin series and truncation in the series representation of each term of the series will have to be used. However, numerical differentiation is avoided and this is still a better approach than evaluating $\wp$ in its original space where it is defined in terms of derivative operations. The error then only depends on the numerical Fourier transform algorithm being implemented and the step size in the propagation direction. The Fourier error can be reduced if certain conditions of the sampling step sizes are maintained (i.e., the Nyquist criterion).

In the exponential form of the $\alpha_1$ operator the function is applied in the original domains without need of a Maclaurin expansion since the argument operator has a range in the complex numbers.

The specific case of how the exponential $\alpha_2$ operator is applied in a manner that reduces operations, while maintaining accuracy and preserving the exponential nature of the operator will be discussed in section 2.4-2.5.

## 2.4 Application of the $\alpha_2$ Operator and its Representation in the Relevant Domains used in its Application

The complication factors in, in how the $\alpha_2$ operator is applied. There are several ways in which this operator can be implemented. The operator can be implemented in the original domain and then truncated after some terms in the Maclaurin series expansion with respect to $\alpha_2$ of the corresponding exponential. This yields time derivative operator terms that act on $u$. Less accurate finite difference methods would be required to evaluate the derivative terms acting on $u$ after the expansion. However, as mentioned in the previous section for greater accuracy and reduced computational requirements, an equivalent version of the $\alpha_2$ operator that has no Maclaurin series expansion required should be used.

The exponential expression of $\alpha_2$ is simplified by first expanding the exponential operator in its Maclaurin series w.r.t to the propagation coordinate (where the step size $\frac{1}{2}\Delta\varsigma$ is included in the constant $c_2$ of $\alpha_2$ for ease of writing):

$$e^{\alpha_2} = 1 + c_2\mathbb{N}\frac{\partial}{\partial\tau} + \frac{1}{2!}(c_2\mathbb{N})^2\frac{\partial^2}{\partial\tau} + \frac{1}{3!}(c_2\mathbb{N})^3\frac{\partial^3}{\partial\tau} \ldots \qquad (21)$$

It is worth noting that the derivative term and $\mathbb{N}$ is assumed to be re-arrangeable in the product terms and power terms in the above expansion. However, $\mathbb{N}$ does not commute with the derivative term so this is not completely true and an approximation. The next section covers this in more detail. For now, this approximation is assumed. $\mathbb{N}$ is calculated with the mean-value approximation outlined in steps 1-7.

Now, it can be seen that ( $R = c_2\mathbb{N}$ ) :

$$e^{\alpha_2(\tau,x,y)}u(\tau,x,y) = \left[1 + R(\tau,x,y)\frac{\partial}{\partial\tau} + \frac{1}{2!}R(\tau,x,y)^2\frac{\partial^2}{\partial\tau} + \frac{1}{3!}R(\tau,x,y)^3\frac{\partial^3}{\partial\tau} \ldots\right]u(\tau,x,y) \qquad (22)$$

The left hand side of Eq. (22) can be obtained by carrying out the following integral operation below[4], where $w'$ is defined as an additional independent variable, and $u(w', x, y)$ as equivalent to the distribution given by the Fourier transform of $u$ to the inverse variable domain of $\tau$ (ie, $u(w', x, y) = u(w, x, y)$ )[5]:

$$\int_{-\infty}^{\infty} \left[ 1 + R(\tau, x, y)(-iw') + \frac{1}{2!}R(\tau, \ldots)^2(-iw')^2 \right.$$
$$\left. + \frac{1}{3!}R(\tau, \ldots)^3(-iw')^3 \ldots \right] u(w', x, y) e^{-iw'\tau} dw' \qquad (23)$$
$$= \left[ 1 + R(\tau, \ldots)\frac{\partial}{\partial \tau} + \frac{1}{2!}R(\tau, \ldots)^2\frac{\partial^2}{\partial \tau} + \frac{1}{3!}R(\tau, \ldots)^3\frac{\partial^3}{\partial \tau} \ldots \right] u(\tau, x, y)$$
$$= e^{\alpha_2(\tau, x, y)} u(\tau, x, y)$$

Therefore,

$$e^{\alpha_2(\tau, \ldots)} u(\tau, \ldots) = \int_{-\infty}^{\infty} \left[ 1 + R(\tau, \ldots)(-iw') + \frac{1}{2!}R(\tau, \ldots)^2(-iw')^2 \right.$$
$$\left. + \frac{1}{3!}R(\tau, \ldots)^3(-iw')^3 \ldots \right] u(w', \ldots) e^{-iw'\tau} dw' \qquad (24)$$

Where, $R$ is unaffected by the Fourier integral since the integral is only over the independent variable $w'$. As well, since $u$ is L2-normed each term in the Fourier integral is bound with the over-arching assumption discussed in section 2.1 that there is an integrable bound solution in L2 to Eq. (1) with the initial input $u$. The overall summation of Fourier integrals in the right hand side of Eq. (24) yields an integrable bound solution in L2.

Eq. (23) can be simplified further by using the Maclaurin series identity to the series within the Fourier integral in Eq. (23), yielding:

$$e^{-c_2 Niw'} = \left[ 1 + R(\tau, \ldots)(-iw') + \frac{1}{2!}R(\tau, \ldots)^2(-iw')^2 + \frac{1}{3!}R(\tau, \ldots)^3(-iw')^3 \ldots \right] \qquad (25)$$

Where, the definition of $R$ on the left hand side of Eq.(25) is used. The argument of the exponential now consists of a function defined on 4 independent variables: $w', \tau, x, y$. From the above, the following relation is obtained through substituting Eq.(25) into Eq. (24):

$$e^{\alpha_2(\tau, \ldots)} u(\tau, \ldots) = f_{w'}\{e^{-c_2 Niw'} u(w', \ldots)\}|_\tau \qquad (26)$$

$f_{w'}$ means the Fourier transform over $w'$ at $\tau$. This identity is used which gets past truncation errors and lets $e^{\alpha_2}$ be applied to $u$ in terms of Fourier transform integrals.

$\overline{\alpha_2}$ is defined as the name of the operator $-c_2 Niw'$ when the above derived substitution identity is used.

---

[4] By the use of Fourier identities for derivatives.
[5] The Fourier constants in front of the integral are omitted for clarity.

The semantics of the $\alpha_2$ operator will become clear once the general steps that are applied with all of these operators in their respective spaces are stated. Let Z be the series of numerical operations:

$$Z = iFFT_{k_x,k_y,w \to x,y,\tau} e^{\frac{1}{4}\wp(k_x,k_y,w)\Delta\varsigma} FFT_{x,y,\tau \to k_x,k_y,w}$$

$$iFFT_{w' \to \tau} e^{\frac{1}{2}\overline{\alpha_2(\chi,\psi,w',\tau)}\Delta\varsigma} iFFT_{k_x,k_y,w \to x,y,\tau} e^{\frac{1}{4}\wp(k_x,k_y,w)\Delta\varsigma} FFT_{x,y,\tau \to k_x,k_y,w} \quad (27)$$

The term $FFT_{w' \to \tau} e^{\frac{1}{2}\overline{\alpha_2(x,y,w',\tau)}\Delta\varsigma}$ shown in Eq. (27) applied to $u$ does the following: $FFT_{w' \to \tau} \left\{ e^{\frac{1}{2}\overline{\alpha_2(x,y,w',\tau)}\Delta\varsigma} u(w',x,y) \right\}$. This means:

1. $u$ is inputted in its spatial-frequency domain representation[6] because its frequency representation is the same as its representation in the $w'$ domain, (i.e., $u(w',x,y) = u(w,x,y)$).
2. The exponent is multiplied into $u(w',c,y)$ across $w'$ at a value of $x,y,\tau$.
3. The value of the inverse Fourier transform on the $w'$ domain of the new function (created in 2.) only at the value of $\tau$ used in 2. is taken. This is the value of the updated $u$ at $x,y,\tau$ coming out of the exponential $\alpha_2$ operator step.
4. The process is repeated for all $x,y,\tau$. At the input of this step $u(x,y,w)$ is sent and after this step an updated $u(x,y,\tau)$ is found[7].

Eq. (27) yields:

$$u(x,y,\tau,\varsigma') = Ze^{\alpha_1 \Delta\varsigma} Zu(x,y,\tau,\varsigma' - \Delta\varsigma) \quad (28)$$

This is done iteratively overall all steps in the propagation coordinate, $\varsigma$.

One can see from the above treatment the power of this exponential operator theory. At the heart of traditional FD methods or Runge-Kutta methods the differential operators are replaced, with a numerical difference scheme. However, this method yields a step size dependent error and is inherently a computational approximation to the derivative terms. However, by using the Fourier representation of derivative terms, there is no violation in the nature of the derivative term, the derivative operator term is simply being replaced by its equivalent algebraic integral representation. There is no step size dependent error in this sense and provided that the Fourier transform can be represented by the FFT algorithm accurately, i.e. if the Nyquist criterion is met, there is no other over all errors in computing these derivatives.

Applying an FD or Runge-Kutta method to solve the differential equations in steps 1-7 also bears more error in the integration than the preceding method. For example, taking step 1, a FD method would look like:

---

[6] $iFFT_{k_x,k_y,w \to x,y,w}$ refers to a two dimensional inverse Fourier transform only over the momentum coordinates.

[7] A 3-D function (representing $u$) is inputted into this exponential $\alpha_2$ operator step. The exponential $\alpha_2$ operator is a 4-D function. During the application of the step, the new function found in 2. Is 4-D. After step 3 is applied over all $x,y,\tau$ the function is a 3-D function.

$$u_{\wp}(\varsigma_{k+\frac{1}{4}}) = u_{\wp}(\varsigma_k) + \frac{\partial u_{\wp}(\varsigma)}{\partial \varsigma}\bigg|\varsigma_k(\frac{1}{4}\Delta\varsigma) = u_{\wp}(\varsigma_k) + \wp u_{\wp}(\varsigma_k)(\frac{1}{4}\Delta\varsigma) \quad (29)$$

Where, the operators can be carried out in the Fourier treatment or with Taylor series expansions (specifically those with differential operators). The more pronounced error in the integration arises because the mean value theorem would still have to be applied for these methods as well, and these methods turn out to be a computational approximation of the analytic solution under the mean value approximation. However, the exponential method outlined in this chapter is the true analytic solution of the integration under the mean value approximation. Runge-Kutta methods bare the same form albeit more developed as Eq. (29). Therefore, it will always out-compete applying Runge-Kutta methods for the longitudinal propagation. The exponential describes the step without making the extra assumption that the system is discretized, it still respects the continuous nature of the problem while other methods do not make this distinction. Adaptive recursive or implicit methods can overcome the mean-value approximation itself, but these are substantially more complex methods that lie out of the scope of this chapter. The computational costs increase and the stability of such methods may be an issue or hard to evaluate.

In sum, the system acts in an exponential manner to the propagation coordinate, if there is only one operator composed of all coefficient terms in Eq. (1) the equation would integrate to the exponential form of that operator[8]. Consequently, the most accurate, stable and intuitive way to model the system would be in an exponential form.

Before concluding this section, another important point must briefly be discussed. The $\mathbb{Q}$ operator presented in section 2.1 is not the most general form of term that this method and the results of this section (and the next section) can be applied to. The reason for defining $\mathbb{Q}$ in section 2.1 was to simplify the math in the subsequent section. However, it can be shown that $\mathbb{Q}$ can be composed of any series of derivative terms (including mixed derivative terms over $x, y, \tau$ and terms to any derivative order) with constant coefficients. The terms that constitute a distributional coefficient to a derivative that arise when $\mathbb{Q}\mathbb{N}$ is expanded can be grouped as the new $\alpha_2$ operator (other terms that emerge in the expansion of $\mathbb{Q}\mathbb{N}$ will be grouped in the new alpha1 operator) and applied accordingly as how the $\alpha_2$ (or $\alpha_1$) operator is outlined in section 2.3-2.4 and the next section.

## 2.5 Approximation Used in Justifying the Series Expansion of the Exponential $\alpha_2$ Operator and Higher Order Update to the Operator

In section 2.4, the Maclaurin series expansion for the $\alpha_2$ operator was shown using the Mean-value approximation outlined in steps 1-7, for $\mathbb{N}$, as:

$$e^{\alpha_2} = 1 + c_2\mathbb{N}\frac{\partial}{\partial\tau} + \frac{1}{2!}(c_2\mathbb{N})^2\frac{\partial^2}{\partial\tau} + \frac{1}{3!}(c_2\mathbb{N})^3\frac{\partial^3}{\partial\tau} \ldots \quad (30)$$

Where the step size is factored into the constant $c_2$. Strictly speaking however, the series expansion of $\alpha_2$ is:

---

[8] Application of the exponential form through the Maclaurin expansion would be extremely complicated as opposed to the 3rd order method shown here based on dividing the equation into a series of operators. However, arbitrary scaling to higher orders can be achieved by this one global operator and the full Maclaurin series represents the exact solution.

$$e^{\alpha_2} = 1 + c_2\mathbb{N}\frac{\partial}{\partial\tau} + \frac{1}{2!}\left(c_2\mathbb{N}\frac{\partial}{\partial\tau}\right)\left(c_2\mathbb{N}\frac{\partial}{\partial\tau}\right) + \frac{1}{3!}\left(c_2\mathbb{N}\frac{\partial}{\partial\tau}\right)\left(c_2\mathbb{N}\frac{\partial}{\partial\tau}\right)\left(c_2\mathbb{N}\frac{\partial}{\partial\tau}\right)\ldots \quad (31)$$

Since $\mathbb{N}$ is a function of $\tau$ the product of the derivatives and $c_2\mathbb{N}$'s cannot be re-arranged. They do not commute.

However, if the approximation is used that $\mathbb{N}$ is a slow-varying function w.r.t to $u$, its derivative w.r.t $\tau$ can be assumed to be negligible (i.e., zero). In that case both of the above equations are equivalent and Eq. (30) can be used. When this approximation cannot be used the need to derive error estimations for $\alpha_2$ and to update the application of the operator in such a way that the error truncation is at least at the same level as general split-step operator methods ($\mathscr{O}(\Delta z^3)$ or higher) is paramount. This is accomplished as follows:

Eq. (31) can be expanded as follows (here to illustrate the step-size exponent, the step size is factored out of $c_2$: $\Delta z = \frac{1}{2}\Delta\varsigma$):

$$e^{\alpha_2} = \left[1 + \Delta z c_2 \mathbb{N}\frac{\partial}{\partial\tau} + \frac{1}{2!}(\Delta z c_2 \mathbb{N})^2 \frac{\partial^2}{\partial\tau} + \frac{1}{3!}(\Delta z c_2 \mathbb{N})^3 \frac{\partial^3}{\partial\tau}\ldots\right] + \frac{\Delta z^2}{2!}\left[(c_2\mathbb{N})\frac{\partial}{\partial\tau}(c_2\mathbb{N})\frac{\partial}{\partial\tau}\right]$$
$$+ \frac{\Delta z^3}{3!}\left[\left[(c_2\mathbb{N})^2\frac{\partial^2}{\partial\tau}(c_2\mathbb{N}) + (c_2\mathbb{N})\left[\frac{\partial}{\partial\tau}(c_2\mathbb{N})\right]^2\right]\frac{\partial}{\partial\tau} + 3(c_2\mathbb{N})^2\frac{\partial}{\partial\tau}(c_2\mathbb{N})\frac{\partial^2}{\partial\tau}\right] \quad (32)$$

The term is corrected up to $\mathscr{O}(\Delta z^4)$. If only $\mathscr{O}(\Delta z^3)$ is desired, then the $\Delta z^3$ term on the RHS of Eq. (32) can be omitted. As described in section 2 since the first term of Eq. (32) is applied in $(x, y, \tau, w')$ space, all corrections are applied in that space as well. Derivatives of $c_2\mathbb{N}$ are numerically evaluated using the Fourier derivative identity. In the appropriate space this yields:

$$e^{\alpha_2(\tau,\ldots)}u(\tau,\ldots) =$$
$$= f_{w'}\left\{\left[e^{-c_2\mathbb{N}iw'} + \frac{\Delta z^2}{2!}\left[(c_2\mathbb{N})\frac{\partial}{\partial\tau}(c_2\mathbb{N})(-iw')\right]\right.\right.$$
$$+ \frac{\Delta z^3}{3!}\left[\left[(c_2\mathbb{N})^2\frac{\partial^2}{\partial\tau}(c_2\mathbb{N}) + (c_2\mathbb{N})\left[\frac{\partial}{\partial\tau}(c_2\mathbb{N})\right]^2\right](-iw')\right.$$
$$\left.\left.+ 3(c_2\mathbb{N})^2\frac{\partial}{\partial\tau}(c_2\mathbb{N})(-iw')^2\right]\right]u(w',\ldots)\right\}|_\tau \quad (33)$$

The $e^{-c_2\mathbb{N}iw'}$ is the series convergence of the first term of Eq. (32) in the $(x, y, \tau, w')$ space and was derived in detail in section 2. The identity from this derivation is used here. Application of the exponential $\alpha_2$ operator in this way can yield a $\mathscr{O}(\Delta z^4)$ error truncation. If additional error truncation is necessary, then additional terms of the expansion of Eq. (31) can be obtained and the summation in Eq. (33) updated accordingly.

## 3 Nyquist Criterion and Pseudo-Spectral Criterion for Spatial Grid Sizes

The Fourier nature of the above derived methodology allows the direct application of the Nyquist sampling criterion. This provides a quantitative approach to reduce step-size error and to adapt grid sizes to counter under-sampling. An adaptive grid size technique was developed for the computational implementation of this method and will be covered in the subsequent white light simulation follow up publication to this paper.

## 4 Conclusions and Extensions

A novel and fast 3+1D numerical technique based on the Strang split-step exponential Fourier method has been developed in this paper. This paper also describes the mathematical techniques of the simulation in [2]. The novel numerical technique herein derived, can solve a wider class of non-linear PDEs (the general form of which is shown in section 2.1), which includes generalized derivative NLSE type equations not accessible with traditional spectral split-step approaches. These generalized derivative NLSE equations are important in optics, notably white-light generation in bulk material. This new approach does not intrinsically have the same numerical errors that are present in Runge-Kutta or finite difference type methods for partial differential equations.

The new method relies on an operator that acts on combined domains (for example, acts both on the time and frequency domain and couples both together). It also extends the linear operator to a more general series case with mixed derivative terms. The approach was shown to maintain the same level of accuracy (to $3^{rd}$ order accuracy) as other Strang spectral split-step methods. A formulism and a method of proof was developed in this paper (extensively in the appendices) that can be used to extend spectral split-step methods further in the breadth of allowed terms and order of accuracy for the corresponding non-linear PDEs.

This paper offers new mathematical tools to simulate non-linear PDEs. The application of these tools to specific non-linear PDEs and computer performance comparisons to other numerical methods will be the focus of upcoming publications and important future steps.

In addition to the analyses presented in section 2, Appendix 1 covers the case where there are convolution terms present in the non-linear partial differential equation. These type of terms arise in physics, for example, when modelling Raman effects.

## Competing Interests

The author does not have any competing interests.

## Acknowledgements

The author would like to thank Prof. Dr. R.J. Dwayne Miller for funding and the opportunity to pursue this project.

## Funding

Funding for this work was supported through the Max Planck Society through the Max Planck Institute for the Structure and Dynamics of Matter.

## Appendix 0: Essential Proofs of Section 2

### 0.0 Fourier convention used

The Fourier convention used is:

$$F(w) = \int_{-\infty}^{\infty} f(t)e^{iwt}dt$$

$$f(t) = \frac{1}{2\pi}\int_{-\infty}^{\infty} F(w)e^{-iwt}dw$$

Where w represents angular frequency and t represents the time domain. The coefficient $\frac{1}{2\pi}$ is omitted in equations and derivations (for clarity) but assumed in this chapter. The convention to represent this integral in quick form is:

$$f(t) = f_w(F(w))|_t$$

## 0.1 Proof: Fourier Representation of series of derivative operators

*This is the proof that any series expansion of derivative operators multiplied into u is equal to the inverse Fourier transform of substituting the derivatives with the representation of the derivatives in Fourier space into the series expansion and multiplying into the representation of u in the Fourier space.*

*In this section the following result will be proved:*

*If*

$$Q = \sum_{n=0}^{\infty}\sum_{j=0}^{\infty} c_{nj}\nabla^n \frac{\partial^j}{\partial^j \tau}$$

$$Q(k_x, k_y, w) = \sum_{n=0}^{\infty}\sum_{j=0}^{\infty} c_{nj}\left((-ik_x)^n + (-ik_y)^n\right)(-iw)^j$$

*Then*

$$Qu(\tau, x, y) = f_{k_x, k_y, w}\left(\left[\sum_{n=0}^{\infty}\sum_{j=0}^{\infty} c_{nj}\left((-ik_x)^n + (-ik_y)^n\right)(-iw)^j\right]u(w, k_x, k_y)\right)\Big|_{\tau, x, y}$$

$$= f_{k_x, k_y, w}\left(Q(k_x, k_y, w)u(w, k_x, k_y)\right)|_{\tau, x, y}$$

If the operator can be represented in the form:

$$Q = \sum_{n=0}^{\infty}\sum_{j=0}^{\infty} c_{nj}\nabla^n \frac{\partial^j}{\partial^j \tau} \tag{34}$$

($c_{nj}$ are constants, the del operator acts over $x, y$ spatial coordinates, or normalized spatial coordinates, $\tau$ is the time coordinate, or normalized time coordinate) which is true for the ℘ operator, then

$$Qu(\tau, x, y) = \sum_{n=0}^{\infty} \sum_{j=0}^{\infty} c_{nj} \nabla^n \left[ \frac{\partial^j u(\tau, x, y)}{\partial^j \tau} \right] \quad (35)$$

The operators $\nabla^n$, $\frac{\partial^j}{\partial^j \tau}$ commute and the order of how they are applied on $u$ does not matter. This follows from the commutation of partial derivatives over independent variables.

Eq.(35) is Fourier transformed in the $w$ domain. Using the Fourier identity for derivative terms to any order, the following is obtained:

$$Qu(\tau, x, y) = \sum_{n=0}^{\infty} \sum_{j=0}^{\infty} c_{nj} \nabla^n \big[ f_w((-iw)^j u(w, x, y))|_\tau \big] \quad (36)$$

Since the integral does not act over $x, y$, $\nabla^n$ can be factored out of the integral. The property of summation of integrals was used to place the Fourier integral within the summation.

Now, to evaluate the del operators over $x, y$ the following 2D Fourier integral can be used:

$$Qu = \sum_{n=0}^{\infty} \sum_{j=0}^{\infty} c_{nj} f_{k_x, k_y} \left( \big((-ik_x)^n + (-ik_y)^n\big) [f_w((-iw)^j u(w, x, y))|_\tau] \right)|_{x,y} \quad (37)$$

Where the Fourier identity of the del operator is used. The summation property of integrals was used.

The integrals commute because they are evaluated over independent variables. The ordering of how the integrals are evaluated does not matter.

Also, due to the summation property of integrals, the integrals can be factored outside of the summation. This simplifies Eq. (37) to:

$$Qu(\tau, x, y) = f_{k_x, k_y, w} \left( \left[ \sum_{n=0}^{\infty} \sum_{j=0}^{\infty} c_{nj} \big((-ik_x)^n + (-ik_y)^n\big) (-iw)^j \right] u(w, k_x, k_y) \right)|_{\tau, x, y} \quad (38)$$

Since the summation of Eq. (38) is simply the Fourier representation of Q,

$$Qu(\tau, x, y) = f_{k_x, k_y, w} \left( Q(k_x, k_y, w) u(w, k_x, k_y) \right)|_{\tau, x, y} \quad (39)$$

Concluding the proof. ∎

$\wp$ is of the same form as Q above, and these proofs apply to $\wp$.

## 0.2 Corollary

*A series expandable function of derivative operators multiplied into $u$, in the original space, equals to the inverse Fourier transform of the same function where the derivative arguments are replaced by the representation of the derivatives in the Fourier space multiplied by the representation of $u$ in the Fourier space.*

The double summation in Eq. (38) mathematically matches the same expansion as Eq. (34) with variable labels replaced. Therefore, the function described is simply the function of the Q operator over derivative arguments replaced by $(-iw)^h$ for the time derivatives and $(-ik_x)^v + (-ik_y)^v$ for the spatial del operator, $\nabla^v$ used in the functional arguments. However, there are important exceptions and conditions to this corollary:

1) Since, the operator definition $\nabla^n \equiv \frac{\partial^n}{\partial^n x} + \frac{\partial^n}{\partial^n y}$ is used, the series convergence to a function over $\nabla$ using the definitions of functions over operators, may not be the same function described by the series where $\nabla^n$ is replaced by $(-ik_x)^n + (-ik_y)^n$ (this is not the case in [2]). This is because, the definition of functions over operators is to replace the function variable by the operator in the Taylor series of the function. However, this definition does not always preserve the algebra of the series over the operator in a domain variable where its range is that of the complex numbers (in the inverse domain): I.e., when the operator is defined in the original domain as in the beginning of 1).
2) The series in Eq. (38), whose range is now the complex numbers, thus may converge only to the function in a specific domain interval. Therefore, the inverse domain interval must be in the region of convergence of the series to the function.

If 1) is true or 2) is not met, the full series form of Eq. (38) will have to be employed without employing any useful properties of the function, for example, the benefits of using the closed form representation-namely that the function value over the inverse domains can be exactly calculated with a finite amount of operations. Provided 1) is true, the imaginary variable substitution (e.g., $-iw$) does not affect the series convergence into the <u>same</u> function: The series of Eq. (34) and subsequently Eq. (38) are power series and if they converge to a given function within an interval then they are also the Taylor series of the function within the same interval. It can be shown that if a Taylor series converges to a function then if an imaginary variable substitution is done on the Taylor series, it will converge to the same function with the imaginary variable substituted for the function variable (provided intervals of convergence are respected, i.e., condition 2)).

A note: If the operator is defined as a series over derivatives that does not have a closed form functional representation, then the above method derived in section 2 is still valid but with a slight difference. The Maclaurin expansion of the exponential operator will be over this series. This means that the full series will be substituted for the variable in the exponential Maclaurin expansion. By the proof in 0.1 and 0.3 (below) the application of the operator to $u$ is equivalent to taking the inverse Fourier transform of the exponential Maclaurin expansion with the Fourier variables substituted in the series multiplied into the Fourier representation of $u$. While this avoids having to do numerical derivatives, the series terms of the Macluarin expansion will need to be truncated, and the general Maclaurin expansion will need to be truncated, generating numerical error.

## 0.3 Proof of the Exponential Fourier Representation of $\wp$

*In this section the following result will be proved:*

$$e^{\wp(\tau,x,y)}u(\tau,x,y) = f_{k_x,k_y,w}\left(e^{\wp(k_x,k_y,w)}u(w,k_x,k_y)\right)|_{\tau,x,y}$$

Expanding the exponential $\wp$ operator w.r.t to the propagation coordinate in a Maclaurin series yields:

$$e^{\wp} = 1 + \wp + \frac{1}{2!}(\wp)^2 + \frac{1}{3!}(\wp)^3 \ldots \tag{40}$$

The step size is factored into $\wp$, for ease of writing. In the above multiplicative ordering (or additive ordering) does not have to be defined since derivative operators commute like numeric variables.

From the above, an equivalent application of the exponential operator goes as:

$$e^{\wp(\tau,\ldots)}u(\tau,\ldots) = [1 + \wp + \frac{1}{2!}(\wp)^2 + \frac{1}{3!}(\wp)^3 \ldots]u(\tau,\ldots) \tag{41}$$

$\wp$ is equivalent to a global summation of derivative terms with constant coefficients as indicated in Eq. (2).

$(\wp(\tau,\ldots))^n$ in expanded form is equivalent to a summation of products of powers of spatial derivative and temporal derivative terms with constant coefficients. Thus, as shown in 0.1, the Fourier identity for derivative terms can be used. Since the now complex numbered variables also commute the expanded form is then factorized in this Fourier space, in such a way that the following is obtained[9]:

$$(\wp(\tau,\ldots))^n u(\tau,\ldots) = f_{k_x,k_y,w}[(\wp(k_x,k_y,w))^n u(k_x,k_y,w)]|_{\tau,x,y} \tag{42}$$

From the identity of Eq. (42), for each term of the RHS of Eq. (41) (which is of the same form as Eq. (42) with different $n$ power values between them), the following is true:

$$\text{term}(\wp(\tau,x,y))u(\tau,x,y) = f_{k_x,k_y,w}[\text{term}(\wp(k_x,k_y,w))u(k_x,k_y,w)]|_{x,y,\tau}$$

Therefore, Eq. (41) is equivalent to:

---

[9] The replacement of the derivative operators with the numerical variables does not change the expression. Nor does how multiplication work change (due to the same commutation between derivative operators and numerical variables). For this point it can be viewed as just replacing symbols and thus, a factorization back into the form $\wp^n$ where $\wp$ is now over the symbol set $k_x, k_y, w$ is trivially possible.

$$e^{\wp(\tau,x,y)}u(\tau,x,y)$$
$$= f_{k_x,k_y,w}\left[\left[1 + \wp(k_x,k_y,w) + \frac{1}{2!}\big(\wp(k_x,k_y,w)\big)^2\right.\right.$$
$$\left.\left.+ \frac{1}{3!}\big(\wp(k_x,k_y,w)\big)^3 \ldots\right]u(k_x,k_y,w)\right]\Big|_{x,y,\tau}))) \tag{43}$$

Using the fact that a summation of integrals over the same domains is equivalent to the integral of the summation over the domain.

The summation on the RHS is the Macluarin series of $e^{\wp(k_x,k_y,w)}$. This yields:

$$e^{\wp(\tau,\ldots)}u(\tau,\ldots) = f_{k_x,k_y,w}\left[e^{\wp(k_x,k_y,w)}u(k_x,k_y,w)\right]\Big|_{x,y,\tau}) \tag{44}$$

Following the discussion in 0.2 $\wp(k_x, k_y, w)$ can be described as a function in this momenta frequency space.

This concludes the proof of the identity. ∎

Also, if the closed form functional representation of $\wp$ in the Fourier space cannot be obtained due to the series not converging to a closed form function, then the full expansion of Eq. (43) will have to be considered, where the series representation of $\wp(k_x, k_y, w)$ is used. Truncating terms both in the series representation of each term based on $\wp(k_x, k_y, w)$ in the Maclaurin series and in the Maclaurin series on the RHS of Eq. (43) will need to be done.

### 0.4 Proof of solution of differential equations in steps 1-7

It can be seen directly by differentiating the Maclaurin series of the exponential operator multiplied into $u$. Since the operator itself does not rely on the propagation coordinate (for example, a mean field approximation is used), the analytic integration of steps 1-7 with an arbitrary operator including convolution operators always yields the exponential operator form as a solution. See appendix 1 for more rigorous details and proofs.

## Appendix 1: Extension of the Model: Raman Terms (Convolution Terms) and Proof of Steps 1-7

### 1.1 Operator for Raman Convolution Terms

In this section, the exponential operator representation in the domain and inverse domains of convolution terms which can arise in optics such as Raman terms will be derived. Also, since the derivation mirrors the proof of the solutions to the differential equations in steps 1-7, while deriving the convolution exponential operator, the proof that can be applied for all differential equations in steps 1-7 of section 2.3 will be showed.

Convolution terms such as:

$$\int_{-\infty}^{\infty} f(t')u(t'-t)\mathrm{d}t' \equiv f(t) \circ u(t) \tag{45}$$

Can still be analytically integrated like in steps 1-7. If $f$ depends on the propagation coordinate (in this appendix denoted as $z$) then the mean field approximation can be used in the same manner as demonstrated in the description of the steps. The rest of the arguments of both functions above $(x, y)$ are not shown for the purposes of clarity. The Maclaurin series expansion of the exponential operator w.r.t the propagation coordinate of Eq. (45) is:

$$e^{\Delta z f(t) \circ} = 1 + \Delta z f(t) \circ + \frac{1}{2!}(\Delta z f(t) \circ)^2 + \frac{1}{3!}(\Delta z f(t) \circ)^3 \ldots \tag{46}$$

It can be verified that this is a solution to the differential equation shown in steps 1-7, as follows:

If $u = e^{zf(t)\circ}g(t)$, and using the identity of Eq. (46) in the differentiation (now with the free variable z), the following is obtained:

$$\frac{\partial e^{zf(t)\circ}g(t)}{\partial z} = [f(t) \circ + \frac{1}{1!}(zf(t)\circ)f(t) \circ + \frac{1}{2!}(zf(t)\circ)^2 f(t) \circ \ldots \frac{1}{3!}(zf(t)\circ)^3 f(t) \circ]g(t) \tag{47}$$

The above was obtained by differentiating Eq. (46) w.r.t to the propagation coordinate variable. Since $f(t)$ is assumed to be only a function of $t$, it can be treated as a constant in the differentiation. Eq. (47) uses the fact that differentiating a term w.r.t the propagation coordinate yields:

$$\frac{\partial \frac{1}{n!}(zf(t)\circ)^n}{\partial z} = \frac{1}{(n-1)!}(z)^{n-1}(f(t)\circ)^n = \frac{1}{(n-1)!}(zf(t)\circ)^{n-1}f(t) \circ \tag{48}$$

Using the normal chain rule of differentiation. Since, convolutions of functions multiplicatively commute by property of convolution [16], the $f(t) \circ$ power terms above can be factored in any way.

$$\frac{1}{(n-1)!}(zf(t)\circ)^{n-1}f(t) \circ = f(t) \circ \frac{1}{(n-1)!}(zf(t)\circ)^{n-1} \tag{49}$$

Using the factorization in Eq. (49), Eq. (47) becomes:

$$\frac{\partial e^{zf(t)\circ}g(t)}{\partial z} = f(t)\circ\left[1+\frac{1}{1!}(zf(t)\circ)+\frac{1}{2!}(zf(t)\circ)^2\ldots\frac{1}{3!}(zf(t)\circ)^3\right]g(t)$$
$$= f(t)\circ(e^{zf(t)\circ}g(t)) = f(t)\circ u(t,z) \tag{50}$$

Where, the last equality follows from the definition of u. If $u = e^{zf(t)\circ}g(t), \frac{\partial u}{\partial z} = f(t)\circ u(t,z)$ which is the differential equation that needs to be satisfied for the convolution term.

At $z = 0, u = u(t,0)$, substituting $z = 0$ above yields:

$$g(t) = u(t,0) \tag{51}$$

Giving, for a general coordinate step the integrated solution at $z = z_0 + \Delta z$ (step-size, $\Delta z$, is weighted by a symmetrisation constant from the global step-size as shown in steps 1-7):

$$u(t,(n+1)\Delta z) = e^{\Delta zf(t)\circ}u(t,n\Delta z) \tag{52}$$

Where $n$ is the integer slice number.

Therefore, the exponential operator for the convolution term is verified to be $e^{\Delta zf(t)\circ}$. If $f(t)$ is truly independent from z, then the above is an exact solution. If $f(t)$ is not independent of $z$ then the above proof is valid with the use of the mean-value approximation. $u = e^{zf(t,z_0)\circ}g(t)$ would have to be used, where $f(t,z_0)$ is treated as independent from z. Using the expansion of the exponential operator w.r.t z about $z_0$, the following is obtained:

$$e^{zf(t,z_0)\circ} = 1 + (z-z_0)f(t,z_0) + \frac{1}{2!}(z-z_0)^2 f(t,z_0)^2 + \frac{1}{3!}(z-z_0)^3 f(t,z_0)^3 \ldots \tag{53}$$

The rest follows the same as the proof above since $f(t,z_0)$ is treated as independent from $z$.

**1.2 Proof of Steps 1-7 For Operators of Section 2**

The above derivation is the same for the $\wp$ operator (given the discussion in appendix 0.3 that ordering of operations in $\wp$ do not need to be considered and $\wp$ is independent of the the propagation coordinate i.e., $\wp$ does not act on z or is a function of z). This Means that $\wp$ can be symbolically substituted for $f(t) \circ$ symbol in the proof above and the $\wp$ differential equation is proved.

The proof above can also be used for $\alpha_1$ operator if the mean-value approximation is used for $\alpha_1$. Because, the power terms of this operator in the Maclaurin exponential series can be multiplicatively factored in any way and are not functions of z.

Just as for the $\alpha_1$ operator, the mean-value approximation is used for the $\alpha_2$ validating the above proof as applied to $\alpha_2$. It can be shown that the ordering of the multiplication of the Maclaurin exponential power series terms can be arbitrary. To illustrate this statement, an example, using a term in the

operator's exponential Maclaurin series taken from Eq. (31), is presented below:

$$\frac{1}{3!}\left(c_2 \mathbb{N} \frac{\partial}{\partial \tau}\right)\left[\left(c_2 \mathbb{N} \frac{\partial}{\partial \tau}\right)\left[\left(c_2 \mathbb{N} \frac{\partial}{\partial \tau}\right) g(\tau)\right]\right] = \frac{1}{3!}\left[\left(c_2 \mathbb{N} \frac{\partial}{\partial \tau}\right)\left(c_2 \mathbb{N} \frac{\partial}{\partial \tau}\right)\right]\left(c_2 \mathbb{N} \frac{\partial}{\partial \tau}\right) g(\tau) \tag{54}$$

$g(\tau)$ is the initial condition function used in the proof of Appendix 1.1 (i.e., it is $u$ coming from the previous step in steps 1-7). The LHS of Eq. (54) means that the term closest to $g$ is applied on $g$ then the middle term is applied and finally the outer term (next to the fraction). The RHS means that the two terms in the bracket are multiplicatively expanded first, the term closest to $g$ is applied on $g$ then the expanded term is applied. The example shown in Eq. (54) is meant to indicate that the multiplicative ordering of terms in the Macluarin series expansion of the exponential $\alpha_2$ does not matter. Thus, it is possible to use the multiplicative ordering defined as what is shown on the LHS of Eq. (54) for each power term of the series without consequence. If this ordering is used the factorization shown in the proof of Appendix 1.1 in Eq. (50) can still be done (with the closest $c_2 \mathbb{N} \frac{\partial}{\partial \tau}$ term to the fraction coefficient of each power term being factored out to the left) since it is equivalent to the application of the power terms with the LHS ordering in Eq. (54). Since this step can be realized, the exponential $\alpha_2$ operator satisfies the general proof outlined in section 1.1 and therefore, solves the $\alpha_2$ differential equation shown in the steps of section 2.

### 1.3 Convolution Term Frequency Representation

Since,

$$e^{f(t)\circ}u = [1 - f(t) \circ + \frac{1}{2!}(f(t) \circ)^2 - \frac{1}{3!}(f(t) \circ)^3 \ldots ]u \tag{55}$$

Therefore, the following is true:

$$e^{f(t)\circ}u = f_W\left[1 - f(w) + \frac{1}{2!}(f(w))^2 - \frac{1}{3!}(f(w))^3 \ldots \right]u(w)|_t \tag{56}$$

For clarity, the symmetrisation specific step-size increment is factored into $f$. The convolution raised to a power is defined as iterative convolutions of function $f$. As discussed, the multiplicative ordering of these convolutions can be arbitrary. Convolution in the time domain is equivalent to the Fourier transform of the product of the Fourier frequency representation of functions in the convolution. For iterative convolutions, it can be shown as a Fourier identity that this is equivalent, in the time domain, to the Fourier transform of the product of the Fourier frequency representation of all functions in the convolution [16].

The above equations are equivalent to:

$$e^{f(t)\circ}u = f_W[e^{f(w)}]u(w)|_t \tag{57}$$

The operator in the frequency space is simply: $e^{f(w)}$ and can be applied in the manner shown in section 2.

### 1.4 Relevant Note on Terminology

In the above analysis the Maclaurin series is used when the coordinate value $z_0$ can be regarded as zero (the origin) for the local coordinate system of the operator without any loss of generality. When using the mean-value approximation, the Maclaurin series can be used with $z_0$ being regarded as the origin without any loss of generality but the reader must keep in mind that the functions of z in the operators are at $z_0$ in the calculation and for that step can be regarded as constant over the propagation coordinate.